\documentclass{amsart}
\usepackage{amssymb}
\usepackage{amsfonts}

\newtheorem{theorem}{Theorem}
\theoremstyle{plain}

\newtheorem{corollary}{Corollary}

\newtheorem{lemma}{Lemma}

\newtheorem{remark}{Remark}

\numberwithin{equation}{section}
\input{tcilatex}

\begin{document}
\title[Gr\"{u}ss Type Inequalities]{Gr\"{u}ss Type Discrete Inequalities in
Inner Product Spaces, Revisited}
\author{Sever S. Dragomir}
\address{School of Computer Science and Mathematics\\
Victoria University of Technology\\
PO Box 14428, MCMC 8001\\
Victoria, Australia.}
\email{sever@matilda.vu.edu.au}
\urladdr{http://rgmia.vu.edu.au/SSDragomirWeb.html}
\date{22 May, 2003.}

\begin{abstract}
Some sharp inequalities of Gr\"{u}ss type for sequences of vectors in real
or complex inner product spaces are obtained. Applications for Jensen's
inequality for convex functions defined on such spaces are also provided.
\end{abstract}

\maketitle

\section{Introduction}

The following inequality of Gr\"{u}ss type for sequences of vectors in inner
product spaces has been established in \cite{D1}.

\begin{theorem}
\label{t1.a}Let $\left( H;\left\langle \cdot ,\cdot \right\rangle \right) $
be an inner product space over the real or complex number field $\mathbb{K}$ 
$\left( \mathbb{K=R},\mathbb{C}\right) ,$ $\overline{\mathbf{x}}=\left(
x_{1},\dots ,x_{n}\right) \in H^{n}$, $\overline{\mathbf{\alpha }}=\left(
\alpha _{1},\dots ,\alpha _{n}\right) \in \mathbb{K}^{n},$ $\overline{%
\mathbf{p}}=\left( p_{1},\dots ,p_{n}\right) \in \mathbb{R}_{+}^{n}$ with $%
\sum_{i=1}^{n}p_{i}=1.$ If $a,A\in \mathbb{K}$ and $x,X\in H$ are such that%
\begin{multline}
\func{Re}\left[ \left( A-\alpha _{i}\right) \left( \overline{\alpha _{i}}-%
\overline{a}\right) \right] \geq 0\text{ \ and \ }\func{Re}\left\langle
X-x_{i},x_{i}-x\right\rangle \geq 0\ \text{\ }  \label{1.1} \\
\text{for each \ }i\in \left\{ 1,\dots ,n\right\} ,
\end{multline}%
then we have the inequality%
\begin{equation}
0\leq \left\Vert \sum_{i=1}^{n}p_{i}\alpha
_{i}x_{i}-\sum_{i=1}^{n}p_{i}\alpha _{i}\cdot
\sum_{i=1}^{n}p_{i}x_{i}\right\Vert \leq \frac{1}{4}\left\vert
A-a\right\vert \left\Vert X-x\right\Vert .  \label{1.2}
\end{equation}%
The constant $\frac{1}{4}$ is the best possible one in the sense that it
cannot be replaced by a smaller constant.
\end{theorem}

Another result of this type is embodied in the following theorem that has
been obtained in \cite{D2}.

\begin{theorem}
\label{t1.2}Let $H,\mathbb{K}$ be as above $\overline{\mathbf{x}}=\left(
x_{1},\dots ,x_{n}\right) ,$ $\overline{\mathbf{y}}=\left( y_{1},\dots
,y_{n}\right) \in H^{n}$ and $\overline{\mathbf{p}}$ a probability sequence. 
$If$ $x,X,y,Y\in H$ are such that 
\begin{multline}
\func{Re}\left\langle X-x_{i},x_{i}-x\right\rangle \geq 0\text{ \ and \ }%
\func{Re}\left\langle Y-y_{i},y_{i}-y\right\rangle \geq 0\ \text{\ }
\label{1.3} \\
\text{for each \ }i\in \left\{ 1,\dots ,n\right\} ,
\end{multline}%
then we have the inequality%
\begin{equation}
0\leq \left\vert \sum_{i=1}^{n}p_{i}\left\langle x_{i},y_{i}\right\rangle
-\left\langle \sum_{i=1}^{n}p_{i}x_{i},\sum_{i=1}^{n}p_{i}y_{i}\right\rangle
\right\vert \leq \frac{1}{4}\left\Vert X-x\right\Vert \left\Vert
Y-y\right\Vert .  \label{1.4}
\end{equation}%
The constant $\frac{1}{4}$ is best possible in the above sense.
\end{theorem}

On choosing $x_{i}=y_{i}$ $\left( i=1,\dots ,n\right) $ in Theorem \ref{t1.2}%
, one may obtain the following counterpart of Cauchy-Bunyakovsky-Schwarz
inequality%
\begin{equation}
0\leq \sum_{i=1}^{n}p_{i}\left\Vert x_{i}\right\Vert ^{2}-\left\Vert
\sum_{i=1}^{n}p_{i}x_{i}\right\Vert ^{2}\leq \frac{1}{4}\left\Vert
X-x\right\Vert ^{2},  \label{1.5}
\end{equation}%
provided $\overline{\mathbf{x}}$ and $\overline{\mathbf{p}}$ satisfy the
assumptions of Theorem \ref{t1.2}.

In the recent paper \cite{D3}, the author has obtained the following Gr\"{u}%
ss type inequality for forward difference as well.

\begin{theorem}
\label{t1.3}Let $\overline{\mathbf{x}}=\left( x_{1},\dots ,x_{n}\right) ,$ $%
\overline{\mathbf{y}}=\left( y_{1},\dots ,y_{n}\right) \in H^{n}$ and $%
\overline{\mathbf{p}}\in \mathbb{R}_{+}^{n}$ be a probability sequence. Then
one has the inequalities%
\begin{multline}
\left\vert \sum_{i=1}^{n}p_{i}\left\langle x_{i},y_{i}\right\rangle
-\left\langle \sum_{i=1}^{n}p_{i}x_{i},\sum_{i=1}^{n}p_{i}y_{i}\right\rangle
\right\vert   \label{1.6} \\
\leq \left\{ 
\begin{array}{l}
\left[ \sum\limits_{i=1}^{n}i^{2}p_{i}-\left(
\sum\limits_{i=1}^{n}ip_{i}\right) ^{2}\right] \max\limits_{k=1,\dots
,n-1}\left\Vert \Delta x_{k}\right\Vert \max\limits_{k=1,\dots
,n-1}\left\Vert \Delta y_{k}\right\Vert ; \\ 
\\ 
\sum\limits_{1\leq j<i\leq n}p_{i}p_{j}\left( i-j\right) \left(
\sum\limits_{k=1}^{n-1}\left\Vert \Delta x_{k}\right\Vert ^{p}\right) ^{%
\frac{1}{p}}\left( \sum\limits_{k=1}^{n-1}\left\Vert \Delta y_{k}\right\Vert
^{q}\right) ^{\frac{1}{q}} \\ 
\hfill \text{if }p>1,\ \frac{1}{p}+\frac{1}{q}=1; \\ 
\\ 
\dfrac{1}{2}\left[ \sum\limits_{i=1}^{n}p_{i}\left( 1-p_{i}\right) \right]
\sum\limits_{k=1}^{n-1}\left\Vert \Delta x_{k}\right\Vert
\sum\limits_{k=1}^{n-1}\left\Vert \Delta y_{k}\right\Vert .%
\end{array}%
\right. 
\end{multline}%
The constants $1,$ $1$ and $\frac{1}{2}$ in the right hand side of the
inequality (\ref{1.6}) are best in the sense that they cannot be replaced by
smaller constants.
\end{theorem}

If one chooses $p_{i}=\frac{1}{n}$ $\left( i=1,\dots ,n\right) $ in (\ref%
{1.6}), then the following unweighted inequalities would hold:%
\begin{multline}
\left\vert \frac{1}{n}\sum_{i=1}^{n}\left\langle x_{i},y_{i}\right\rangle
-\left\langle \frac{1}{n}\sum_{i=1}^{n}x_{i},\frac{1}{n}\sum_{i=1}^{n}y_{i}%
\right\rangle \right\vert   \label{1.7} \\
\leq \left\{ 
\begin{array}{l}
\dfrac{n^{2}-1}{12}\max\limits_{k=1,\dots ,n-1}\left\Vert \Delta
x_{k}\right\Vert \max\limits_{k=1,\dots ,n-1}\left\Vert \Delta
y_{k}\right\Vert ; \\ 
\\ 
\dfrac{n^{2}-1}{6n}\left( \sum\limits_{k=1}^{n-1}\left\Vert \Delta
x_{k}\right\Vert ^{p}\right) ^{\frac{1}{p}}\left(
\sum\limits_{k=1}^{n-1}\left\Vert \Delta y_{k}\right\Vert ^{q}\right) ^{%
\frac{1}{q}} \\ 
\hfill \text{if }p>1,\ \frac{1}{p}+\frac{1}{q}=1; \\ 
\\ 
\dfrac{n-1}{2n}\sum\limits_{k=1}^{n-1}\left\Vert \Delta x_{k}\right\Vert
\sum\limits_{k=1}^{n-1}\left\Vert \Delta y_{k}\right\Vert .%
\end{array}%
\right. 
\end{multline}%
Here, the constants $\frac{1}{12},$ $\frac{1}{6}$ and $\frac{1}{2}$ are also
best possible in the above sense.

The following counterpart inequality of the Cauchy-Bunyakovsky-Schwarz
inequality for sequences of vectors in inner product spaces holds.

\begin{corollary}
\label{c1.4}With the assumptions in Theorem \ref{t1.3} for $\overline{%
\mathbf{x}}$ and $\overline{\mathbf{p}}$ one has the inequalities%
\begin{align}
0& \leq \sum_{i=1}^{n}p_{i}\left\Vert x_{i}\right\Vert ^{2}-\left\Vert
\sum_{i=1}^{n}p_{i}x_{i}\right\Vert ^{2}  \label{1.8} \\
& \leq \left\{ 
\begin{array}{l}
\left[ \sum\limits_{i=1}^{n}i^{2}p_{i}-\left(
\sum\limits_{i=1}^{n}ip_{i}\right) ^{2}\right] \max\limits_{k=\overline{1,n-1%
}}\left\Vert \Delta x_{k}\right\Vert ^{2}; \\ 
\\ 
\sum\limits_{1\leq j<i\leq n}p_{i}p_{j}\left( i-j\right) \left(
\sum\limits_{k=1}^{n-1}\left\Vert \Delta x_{k}\right\Vert ^{p}\right) ^{%
\frac{1}{p}}\left( \sum\limits_{k=1}^{n-1}\left\Vert \Delta x_{k}\right\Vert
^{q}\right) ^{\frac{1}{q}} \\ 
\hfill \text{if }p>1,\ \frac{1}{p}+\frac{1}{q}=1; \\ 
\\ 
\dfrac{1}{2}\left[ \sum\limits_{i=1}^{n}p_{i}\left( 1-p_{i}\right) \right]
\left( \sum\limits_{k=1}^{n-1}\left\Vert \Delta x_{k}\right\Vert \right)
^{2}.%
\end{array}%
\right.   \notag
\end{align}%
The constants $1,$ $1$ and $\frac{1}{2}$ are best possible in the above
sense.
\end{corollary}

The following particular inequalities that may be deduced from (\ref{1.8})
on choosing the equal weights $p_{i}=\frac{1}{n},$ $i=1,\dots ,n$ are also
of interest%
\begin{align}
0& \leq \frac{1}{n}\sum_{i=1}^{n}\left\Vert x_{i}\right\Vert ^{2}-\left\Vert 
\frac{1}{n}\sum_{i=1}^{n}x_{i}\right\Vert ^{2}  \label{1.9} \\
& \leq \left\{ 
\begin{array}{l}
\dfrac{n^{2}-1}{12}\max\limits_{k=\overline{1,n-1}}\left\Vert \Delta
x_{k}\right\Vert ^{2}; \\ 
\\ 
\dfrac{n^{2}-1}{6n}\left( \sum\limits_{k=1}^{n-1}\left\Vert \Delta
x_{k}\right\Vert ^{p}\right) ^{\frac{1}{p}}\left(
\sum\limits_{k=1}^{n-1}\left\Vert \Delta x_{k}\right\Vert ^{q}\right) ^{%
\frac{1}{q}} \\ 
\hfill \text{if }p>1,\ \frac{1}{p}+\frac{1}{q}=1; \\ 
\\ 
\dfrac{n-1}{2n}\left( \sum\limits_{k=1}^{n-1}\left\Vert \Delta
x_{k}\right\Vert \right) ^{2}.%
\end{array}%
\right.   \notag
\end{align}%
Here the constants $\frac{1}{12},$ $\frac{1}{6}$ and $\frac{1}{2}$ are also
best possible.

It is the main aim of this paper to point out a different class of Gr\"{u}ss
type inequalities for sequences of vectors in inner product spaces and to
apply them for obtaining a reverse of Jenssen's inequality for convex
functions defined on such spaces.

\section{Some Gr\"{u}ss Type Inequalities}

The following lemma holds (see also \cite{D4}).

\begin{lemma}
\label{l2.1}Let $a,x,A$ be vectors in the inner product space $\left(
H;\left\langle \cdot ,\cdot \right\rangle \right) $ over the real or complex
number field $\mathbb{K}$ $\left( \mathbb{K=R},\mathbb{C}\right) $ with $%
a\neq A.$ The following statements are equivalent:

\begin{enumerate}
\item[(i)] $\func{Re}\left\langle A-x,x-a\right\rangle \geq 0;$

\item[(ii)] $\left\Vert x-\frac{a+A}{2}\right\Vert \leq \frac{1}{2}%
\left\Vert A-a\right\Vert .$
\end{enumerate}
\end{lemma}

\begin{proof}
For the sake of completeness, we give a simple proof as follows. 

Let%
\begin{equation*}
I_{1}:=\func{Re}\left\langle A-x,x-a\right\rangle =-\func{Re}\left\langle
A,a\right\rangle -\left\Vert x\right\Vert ^{2}+\func{Re}\left[ \left\langle 
\overline{x,A}\right\rangle +\left\langle x,a\right\rangle \right] 
\end{equation*}%
and%
\begin{align*}
I_{2}& :=\frac{1}{4}\left\Vert A-a\right\Vert ^{2}-\left\Vert x-\frac{a+A}{2}%
\right\Vert ^{2} \\
& =\frac{1}{4}\left( \left\Vert A\right\Vert ^{2}-2\func{Re}\left\langle
A,a\right\rangle +\left\Vert a\right\Vert ^{2}\right) -\left( \left\Vert
x\right\Vert ^{2}-\func{Re}\left\langle x,a+A\right\rangle +\frac{1}{4}%
\left\Vert A-a\right\Vert ^{2}\right)  \\
& =-\func{Re}\left\langle A,a\right\rangle -\left\Vert x\right\Vert ^{2}+%
\func{Re}\left[ \left\langle x,a\right\rangle +\left\langle x,A\right\rangle %
\right] .
\end{align*}%
Since $\func{Re}\left\langle x,A\right\rangle =\left[ \func{Re}\overline{%
\left\langle x,A\right\rangle }\right] ,$ we deduce that $I_{1}=I_{2},$
showing the desired equivalence.
\end{proof}

\begin{remark}
\label{r2.2}If $H=\mathbb{C}$, $\left\Vert \cdot \right\Vert =\left\vert
\cdot \right\vert ,$ then the following sentences are equivalent

\begin{enumerate}
\item[(a)] $\func{Re}\left[ \left( A-x\right) \left( \overline{x}-\overline{a%
}\right) \right] \geq 0;$

\item[(aa)] $\left\vert x-\frac{a+A}{2}\right\vert \leq \frac{1}{2}%
\left\vert A-a\right\vert ,$
\end{enumerate}

where $a,A,x\in \mathbb{C}.$

If $H=\mathbb{R}$, $\left\Vert \cdot \right\Vert =\left\vert \cdot
\right\vert ,$ and $A>a,$ then the following sentences are obviously
equivalent:

\begin{enumerate}
\item[(b)] $a\leq x\leq A;$

\item[(bb)] $\left\vert x-\frac{a+A}{2}\right\vert \leq \frac{1}{2}%
\left\vert A-a\right\vert .$
\end{enumerate}
\end{remark}

The following inequality of Gr\"{u}ss type for sequences of vectors in inner
product spaces holds.

\begin{theorem}
\label{t2.3}Let $\left( H;\left\langle \cdot ,\cdot \right\rangle \right) $
be an inner product over $\mathbb{K}$ $\left( \mathbb{K=C},\mathbb{R}\right) 
$, and $\overline{\mathbf{x}}=\left( x_{1},\dots ,x_{n}\right) ,$ $\overline{%
\mathbf{y}}=\left( y_{1},\dots ,y_{n}\right) \in H^{n},$ $\overline{\mathbf{p%
}}\in \mathbb{R}_{+}^{n}$ with $\sum_{i=1}^{n}p_{i}=1.$ If $x,X\in H$ are
such that%
\begin{equation}
\func{Re}\left\langle X-x_{i},x_{i}-x\right\rangle \geq 0\ \text{\ for each
\ }i\in \left\{ 1,\dots ,n\right\} ,  \label{2.1}
\end{equation}%
or, equivalently,%
\begin{equation}
\left\Vert x_{i}-\frac{x+X}{2}\right\Vert \leq \frac{1}{2}\left\Vert
X-x\right\Vert \ \text{\ for each \ }i\in \left\{ 1,\dots ,n\right\} ,
\label{2.2}
\end{equation}%
then one has the inequality%
\begin{align}
& \left\vert \sum_{i=1}^{n}p_{i}\left\langle x_{i},y_{i}\right\rangle
-\left\langle \sum_{i=1}^{n}p_{i}x_{i},\sum_{i=1}^{n}p_{i}y_{i}\right\rangle
\right\vert  \label{2.3} \\
& \leq \frac{1}{2}\left\Vert X-x\right\Vert \sum_{i=1}^{n}p_{i}\left\Vert
y_{i}-\sum_{j=1}^{n}p_{j}y_{j}\right\Vert  \notag \\
& \leq \frac{1}{2}\left\Vert X-x\right\Vert \left[ \sum_{i=1}^{n}p_{i}\left%
\Vert y_{i}\right\Vert ^{2}-\left\Vert \sum_{i=1}^{n}p_{i}y_{i}\right\Vert
^{2}\right] ^{\frac{1}{2}}.  \notag
\end{align}%
The constant $\frac{1}{2}$ is best possible in the first and second
inequality in the sense that it cannot be replaced by a smaller constant.
\end{theorem}

\begin{proof}
It is easy to see that the following identity holds true%
\begin{equation}
\sum_{i=1}^{n}p_{i}\left\langle x_{i},y_{i}\right\rangle -\left\langle
\sum_{i=1}^{n}p_{i}x_{i},\sum_{i=1}^{n}p_{i}y_{i}\right\rangle
=\sum_{i=1}^{n}p_{i}\left\langle x_{i}-\frac{x+X}{2},y_{i}-%
\sum_{j=1}^{n}p_{j}y_{j}\right\rangle .  \label{2.4}
\end{equation}%
Taking the modulus in (\ref{2.4}) and using the Schwarz inequality in the
inner product space $\left( H;\left\langle \cdot ,\cdot \right\rangle
\right) ,$ we have%
\begin{align*}
\left\vert \sum_{i=1}^{n}p_{i}\left\langle x_{i},y_{i}\right\rangle
-\left\langle \sum_{i=1}^{n}p_{i}x_{i},\sum_{i=1}^{n}p_{i}y_{i}\right\rangle
\right\vert & \leq \sum_{i=1}^{n}p_{i}\left\vert \left\langle x_{i}-\frac{x+X%
}{2},y_{i}-\sum_{j=1}^{n}p_{j}y_{j}\right\rangle \right\vert \\
& \leq \sum_{i=1}^{n}p_{i}\left\Vert x_{i}-\frac{x+X}{2}\right\Vert
\left\Vert y_{i}-\sum_{j=1}^{n}p_{j}y_{j}\right\Vert \\
& \leq \frac{1}{2}\left\Vert X-x\right\Vert \sum_{i=1}^{n}p_{i}\left\Vert
y_{i}-\sum_{j=1}^{n}p_{j}y_{j}\right\Vert ,
\end{align*}%
and the first inequality in (\ref{2.3}) is proved.

Using the Cauchy-Bunyakovsky-Schwarz inequality for positive sequences and
the calculation rules in inner product spaces, we have%
\begin{equation*}
\sum_{i=1}^{n}p_{i}\left\Vert y_{i}-\sum_{j=1}^{n}p_{j}y_{j}\right\Vert \leq 
\left[ \sum_{i=1}^{n}p_{i}\left\Vert
y_{i}-\sum_{j=1}^{n}p_{j}y_{j}\right\Vert ^{2}\right] ^{\frac{1}{2}}
\end{equation*}%
and%
\begin{equation*}
\sum_{i=1}^{n}p_{i}\left\Vert y_{i}-\sum_{j=1}^{n}p_{j}y_{j}\right\Vert
^{2}=\sum_{i=1}^{n}p_{i}\left\Vert y_{i}\right\Vert ^{2}-\left\Vert
\sum_{i=1}^{n}p_{i}y_{i}\right\Vert ^{2}
\end{equation*}%
giving the second part of (\ref{2.3}).

To prove the sharpness of the constant $\frac{1}{2}$ in the first inequality
in (\ref{2.3}), let us assume that, under the assumptions of the theorem,
the inequality holds with a constant $C>0,$ i.e., 
\begin{equation}
\left\vert \sum_{i=1}^{n}p_{i}\left\langle x_{i},y_{i}\right\rangle
-\left\langle \sum_{i=1}^{n}p_{i}x_{i},\sum_{i=1}^{n}p_{i}y_{i}\right\rangle
\right\vert \leq C\left\Vert X-x\right\Vert \sum_{i=1}^{n}p_{i}\left\Vert
y_{i}-\sum_{j=1}^{n}p_{j}y_{j}\right\Vert .  \label{2.4a}
\end{equation}%
Consider $n=2$ and observe that%
\begin{align*}
\sum_{i=1}^{2}p_{i}\left\langle x_{i},y_{i}\right\rangle -\left\langle
\sum_{i=1}^{2}p_{i}x_{i},\sum_{i=1}^{2}p_{i}y_{i}\right\rangle &
=p_{2}p_{1}\left\langle x_{2}-x_{1},y_{2}-y_{1}\right\rangle , \\
\sum_{i=1}^{2}p_{i}\left\Vert y_{i}-\sum_{j=1}^{2}p_{j}y_{j}\right\Vert &
=2p_{2}p_{1}\left\Vert y_{2}-y_{1}\right\Vert
\end{align*}%
and then, by (\ref{2.4a}), we deduce%
\begin{equation}
p_{2}p_{1}\left\vert \left\langle x_{2}-x_{1},y_{2}-y_{1}\right\rangle
\right\vert \leq 2C\left\Vert X-x\right\Vert p_{2}p_{1}\left\Vert
y_{2}-y_{1}\right\Vert .  \label{2.5}
\end{equation}%
If we choose $p_{1},p_{2}>0,$ $y_{2}=x_{2},$ $y_{1}=x_{1}$ and $x_{2}=X,$ $%
x_{1}=x$ with $x\neq X,$ then (\ref{2.2}) holds and from (\ref{2.5}) we
deduce $C\geq \frac{1}{2}.$

The fact that $\frac{1}{2}$ is best possible in the second inequality may be
proven in a similar manner and we omit the details.
\end{proof}

\begin{remark}
\label{r2.4}If $\overline{\mathbf{x}}$ and $\overline{\mathbf{y}}$ satisfy
the assumptions of Theorem \ref{t1.2}, or equivalently%
\begin{equation}
\left\Vert x_{i}-\frac{x+X}{2}\right\Vert \leq \frac{1}{2}\left\Vert
X-x\right\Vert ,\ \ \ \ \left\Vert y_{i}-\frac{y+Y}{2}\right\Vert \leq \frac{%
1}{2}\left\Vert Y-y\right\Vert ,  \label{2.6}
\end{equation}%
for each $i\in \left\{ 1,\dots ,n\right\} ,$ then by Theorem \ref{t2.3} we
may state the following sequence of inequalities improving the Gr\"{u}ss
inequality (\ref{2.4})%
\begin{align}
0& \leq \left\vert \sum_{i=1}^{n}p_{i}\left\langle x_{i},y_{i}\right\rangle
-\left\langle \sum_{i=1}^{n}p_{i}x_{i},\sum_{i=1}^{n}p_{i}y_{i}\right\rangle
\right\vert  \label{2.7} \\
& \leq \frac{1}{2}\left\Vert X-x\right\Vert \sum_{i=1}^{n}p_{i}\left\Vert
y_{i}-\sum_{j=1}^{n}p_{j}y_{j}\right\Vert  \notag \\
& \leq \frac{1}{2}\left\Vert X-x\right\Vert \left(
\sum_{i=1}^{n}p_{i}\left\Vert y_{i}\right\Vert ^{2}-\left\Vert
\sum_{i=1}^{n}p_{i}y_{i}\right\Vert ^{2}\right) ^{\frac{1}{2}}  \notag \\
& \leq \frac{1}{4}\left\Vert X-x\right\Vert \left\Vert Y-y\right\Vert .
\end{align}%
In particular, for $x_{i}=y_{i}$ $\left( i=1,\dots ,n\right) ,$ one has%
\begin{equation}
0\leq \sum_{i=1}^{n}p_{i}\left\Vert x_{i}\right\Vert ^{2}-\left\Vert
\sum_{i=1}^{n}p_{i}x_{i}\right\Vert ^{2}\leq \frac{1}{2}\left\Vert
X-x\right\Vert \sum_{i=1}^{n}p_{i}\left\Vert
x_{i}-\sum_{j=1}^{n}p_{j}x_{j}\right\Vert  \label{2.8}
\end{equation}%
and the constant $\frac{1}{2}$ is best possible.
\end{remark}

The following result is connected to Theorem \ref{t1.a} from Introduction.

\begin{theorem}
\label{t2.5}Let $\left( H;\left\langle \cdot ,\cdot \right\rangle \right) $
and $\mathbb{K}$ be as above and $\overline{\mathbf{x}}=\left( x_{1},\dots
,x_{n}\right) \in H^{n},$ $\overline{\mathbf{\alpha }}=\left( \alpha
_{1},\dots ,\alpha _{n}\right) \in \mathbb{K}^{n}$ and $\overline{\mathbf{p}}
$ a probability vector. If $x,X\in H$ are such that (\ref{2.1}) or,
equivalently, (\ref{2.2}) holds, then we have the inequality%
\begin{align}
0& \leq \left\Vert \sum_{i=1}^{n}p_{i}\alpha
_{i}x_{i}-\sum_{i=1}^{n}p_{i}\alpha _{i}\cdot
\sum_{i=1}^{n}p_{i}x_{i}\right\Vert  \label{2.9} \\
& \leq \frac{1}{2}\left\Vert X-x\right\Vert \sum_{i=1}^{n}p_{i}\left\vert
\alpha _{i}-\sum_{j=1}^{n}p_{j}\alpha _{j}\right\vert  \notag \\
& \leq \frac{1}{2}\left\Vert X-x\right\Vert \left[ \sum_{i=1}^{n}p_{i}\left%
\vert \alpha _{i}\right\vert ^{2}-\left\vert \sum_{i=1}^{n}p_{i}\alpha
_{i}\right\vert ^{2}\right] ^{\frac{1}{2}}.  \notag
\end{align}%
The constant $\frac{1}{2}$ in the first and second inequalities is best
possible in the sense that it cannot be replaced by a smaller constant.
\end{theorem}

\begin{proof}
We start with the following equality that may be easily verified by direct
calculation%
\begin{equation}
\sum_{i=1}^{n}p_{i}\alpha _{i}x_{i}-\sum_{i=1}^{n}p_{i}\alpha _{i}\cdot
\sum_{i=1}^{n}p_{i}x_{i}=\sum_{i=1}^{n}p_{i}\left( \alpha
_{i}-\sum_{j=1}^{n}p_{j}\alpha _{j}\right) \left( x_{i}-\frac{x+X}{2}\right)
.  \label{2.10}
\end{equation}%
If we take the norm in (\ref{2.10}), we deduce%
\begin{align*}
\left\Vert \sum_{i=1}^{n}p_{i}\alpha _{i}x_{i}-\sum_{i=1}^{n}p_{i}\alpha
_{i}\cdot \sum_{i=1}^{n}p_{i}x_{i}\right\Vert & \leq
\sum_{i=1}^{n}p_{i}\left\vert \alpha _{i}-\sum_{j=1}^{n}p_{j}\alpha
_{j}\right\vert \left\Vert x_{i}-\frac{x+X}{2}\right\Vert  \\
& \leq \frac{1}{2}\left\Vert X-x\right\Vert \sum_{i=1}^{n}p_{i}\left\vert
\alpha _{i}-\sum_{j=1}^{n}p_{j}\alpha _{j}\right\vert  \\
& \leq \frac{1}{2}\left\Vert X-x\right\Vert \left( \sum_{i=1}^{n}p_{i}\left(
\alpha _{i}-\sum_{j=1}^{n}p_{j}\alpha _{j}\right) ^{2}\right) ^{\frac{1}{2}}
\\
& =\frac{1}{2}\left\Vert X-x\right\Vert \left( \sum_{i=1}^{n}p_{i}\left\vert
\alpha _{i}\right\vert ^{2}-\left\vert \sum_{i=1}^{n}p_{i}\alpha
_{i}\right\vert ^{2}\right) ^{\frac{1}{2}},
\end{align*}%
proving the inequality (\ref{2.9}).

The fact that the constant $\frac{1}{2}$ is sharp may be proven in a similar
manner to the one embodied in the proof of Theorem \ref{t2.3}. We omit the
details.
\end{proof}

\begin{remark}
\label{r2.6}If $\overline{\mathbf{x}}$ and $\overline{\mathbf{\alpha }}$
satisfy the assumption of Theorem \ref{t1.2}, or, equivalently,%
\begin{equation*}
\left\Vert \alpha _{i}-\frac{a+A}{2}\right\Vert \leq \frac{1}{2}\left\vert
A-a\right\vert ,\ \ \ \ \left\Vert x_{i}-\frac{x+X}{2}\right\Vert \leq \frac{%
1}{2}\left\Vert X-x\right\Vert ,
\end{equation*}%
for each $i\in \left\{ 1,\dots ,n\right\} ,$ then by Theorem \ref{t2.5} we
may state the following sequence of inequalities improving the Gr\"{u}ss
inequality (\ref{1.2}),%
\begin{align}
0& \leq \left\Vert \sum_{i=1}^{n}p_{i}\alpha
_{i}x_{i}-\sum_{i=1}^{n}p_{i}\alpha _{i}\cdot
\sum_{i=1}^{n}p_{i}x_{i}\right\Vert  \label{2.11} \\
& \leq \frac{1}{2}\left\Vert X-x\right\Vert \sum_{i=1}^{n}p_{i}\left\vert
\alpha _{i}-\sum_{j=1}^{n}p_{j}\alpha _{j}\right\vert  \notag \\
& \leq \frac{1}{2}\left\Vert X-x\right\Vert \left(
\sum_{i=1}^{n}p_{i}\left\vert \alpha _{i}\right\vert ^{2}-\left\vert
\sum_{i=1}^{n}p_{i}\alpha _{i}\right\vert ^{2}\right) ^{\frac{1}{2}}  \notag
\\
& \leq \frac{1}{4}\left\vert A-a\right\vert \left\Vert X-x\right\Vert . 
\notag
\end{align}
\end{remark}

\begin{remark}
\label{r2.7}If in (\ref{2.9}) we choose $x_{i}=\alpha _{i}\in \mathbb{C}$
and assume that $\left\vert \alpha _{i}-\frac{a+A}{2}\right\vert \leq \frac{1%
}{2}\left\vert A-a\right\vert ,$ where $a,A\in \mathbb{C}$, then we get the
following interesting inequality for complex numbers%
\begin{align*}
0& \leq \left\vert \sum_{i=1}^{n}p_{i}\alpha _{i}^{2}-\left(
\sum_{i=1}^{n}p_{i}\alpha _{i}\right) ^{2}\right\vert \\
& \leq \frac{1}{2}\left\vert A-a\right\vert \sum_{i=1}^{n}p_{i}\left\vert
\alpha _{i}-\sum_{j=1}^{n}p_{j}\alpha _{j}\right\vert \\
& \leq \frac{1}{2}\left\vert A-a\right\vert \left[ \sum_{i=1}^{n}p_{i}\left%
\vert \alpha _{i}\right\vert ^{2}-\left\vert \sum_{i=1}^{n}p_{i}\alpha
_{i}\right\vert ^{2}\right] ^{\frac{1}{2}}.
\end{align*}
\end{remark}

\section{Applications for Convex Functions}

Let $\left( H;\left\langle \cdot ,\cdot \right\rangle \right) $ be a real
inner product space and $F:H\rightarrow \mathbb{R}$ a Fr\'{e}chet
differentiable convex function on $H.$ If $\triangledown F:H\rightarrow H$
denotes the gradient operator associated to $F,$ then we have the inequality%
\begin{equation}
F\left( x\right) -F\left( y\right) \geq \left\langle \triangledown F\left(
y\right) ,x-y\right\rangle   \label{3.1}
\end{equation}%
for each $x,y\in H.$

The following result holds.

\begin{theorem}
\label{t3.1}Let $F:H\rightarrow \mathbb{R}$ be as above and $z_{i}\in H,$ $%
i\in \left\{ 1,\dots ,n\right\} .$ Suppose that there exists the vectors $%
m,M\in H$ such that either%
\begin{equation}
\left\langle \triangledown F\left( z_{i}\right) -m,M-\triangledown F\left(
z_{i}\right) \right\rangle \geq 0\text{ \ for each \ }i\in \left\{ 1,\dots
,n\right\} ;  \label{3.2}
\end{equation}%
or, equivalently,%
\begin{equation}
\left\Vert \triangledown F\left( z_{i}\right) -\frac{m+M}{2}\right\Vert \leq 
\frac{1}{2}\left\Vert M-m\right\Vert \text{ \ for each \ }i\in \left\{
1,\dots ,n\right\} .  \label{3.3}
\end{equation}%
If $q_{i}\geq 0$ \ $\left( i\in \left\{ 1,\dots ,n\right\} \right) $ with $%
Q_{n}:=\sum_{i=1}^{n}q_{i}>0,$ then we have the following converse of
Jensen's inequality%
\begin{align}
0& \leq \frac{1}{Q_{n}}\sum_{i=1}^{n}q_{i}F\left( z_{i}\right) -F\left( 
\frac{1}{Q_{n}}\sum_{i=1}^{n}q_{i}z_{i}\right)   \label{3.4} \\
& \leq \frac{1}{2}\left\Vert M-m\right\Vert \frac{1}{Q_{n}}%
\sum_{i=1}^{n}q_{i}\left\Vert z_{i}-\frac{1}{Q_{n}}\sum_{j=1}^{n}q_{j}z_{j}%
\right\Vert   \notag \\
& \leq \frac{1}{2}\left\Vert M-m\right\Vert \left[ \frac{1}{Q_{n}}%
\sum_{i=1}^{n}q_{i}\left\Vert z_{i}\right\Vert ^{2}-\left\Vert \frac{1}{Q_{n}%
}\sum_{i=1}^{n}q_{i}z_{i}\right\Vert ^{2}\right] ^{\frac{1}{2}}.  \notag
\end{align}
\end{theorem}

\begin{proof}
We know, see for example \cite[Eq. (4.4)]{D2}, the following counterpart of
Jensen's inequality for Fr\'{e}chet differentiable convex functions%
\begin{align}
0& \leq \frac{1}{Q_{n}}\sum_{i=1}^{n}q_{i}F\left( z_{i}\right) -F\left( 
\frac{1}{Q_{n}}\sum_{i=1}^{n}q_{i}z_{i}\right)   \label{3.5} \\
& \leq \frac{1}{Q_{n}}\sum_{i=1}^{n}q_{i}\left\langle \triangledown F\left(
z_{i}\right) ,z_{i}\right\rangle -\left\langle \frac{1}{Q_{n}}%
\sum_{i=1}^{n}q_{i}\triangledown F\left( z_{i}\right) ,\frac{1}{Q_{n}}%
\sum_{i=1}^{n}q_{i}z_{i}\right\rangle   \notag
\end{align}%
holds.

Now, if we use Theorem \ref{t2.3} for the choices $x_{i}=\nabla F\left(
z_{i}\right) ,$ $y_{i}=z_{i}$ and $p_{i}=\frac{1}{Q_{n}}q_{i},$ $i\in
\left\{ 1,\dots ,n\right\} ,$ then we can state the inequality%
\begin{align}
& \frac{1}{Q_{n}}\sum_{i=1}^{n}q_{i}\left\langle \triangledown F\left(
z_{i}\right) ,z_{i}\right\rangle -\left\langle \frac{1}{Q_{n}}%
\sum_{i=1}^{n}q_{i}\triangledown F\left( z_{i}\right) ,\frac{1}{Q_{n}}%
\sum_{i=1}^{n}q_{i}z_{i}\right\rangle  \label{3.6} \\
& \leq \frac{1}{2}\left\Vert M-m\right\Vert \frac{1}{Q_{n}}%
\sum_{i=1}^{n}q_{i}\left\Vert z_{i}-\frac{1}{Q_{n}}\sum_{j=1}^{n}q_{j}z_{j}%
\right\Vert  \notag \\
& \leq \frac{1}{2}\left\Vert M-m\right\Vert \left[ \frac{1}{Q_{n}}%
\sum_{i=1}^{n}q_{i}\left\Vert z_{i}\right\Vert ^{2}-\left\Vert \frac{1}{Q_{n}%
}\sum_{i=1}^{n}q_{i}z_{i}\right\Vert ^{2}\right] ^{\frac{1}{2}}.  \notag
\end{align}%
Utilizing (\ref{3.5}) and (\ref{3.6}), we deduce the desired result (\ref%
{3.4}).
\end{proof}

If more information is available about the vector sequence $\overline{%
\mathbf{z}}=\left( z_{1},\dots ,z_{n}\right) \in H^{n},$ then we may state
the following corollary.

\begin{corollary}
\label{c3.2}With the assumptions in Theorem \ref{t3.1} and if there exists
the vectors $z,Z\in H$ such that either%
\begin{equation}
\left\langle z_{i}-z,Z-z_{i}\right\rangle \geq 0\text{ \ for each \ }i\in
\left\{ 1,\dots ,n\right\} ;  \label{3.7}
\end{equation}%
or, equivalently%
\begin{equation}
\left\Vert z_{i}-\frac{z+Z}{2}\right\Vert \leq \frac{1}{2}\left\Vert
Z-z\right\Vert \text{ \ for each \ }i\in \left\{ 1,\dots ,n\right\} ,
\label{3.8}
\end{equation}%
then we have the inequality%
\begin{align}
0& \leq \frac{1}{Q_{n}}\sum_{i=1}^{n}q_{i}F\left( z_{i}\right) -F\left( 
\frac{1}{Q_{n}}\sum_{i=1}^{n}q_{i}z_{i}\right)   \label{3.9} \\
& \leq \frac{1}{2}\left\Vert M-m\right\Vert \frac{1}{Q_{n}}%
\sum_{i=1}^{n}q_{i}\left\Vert z_{i}-\frac{1}{Q_{n}}\sum_{j=1}^{n}q_{j}z_{j}%
\right\Vert   \notag \\
& \leq \frac{1}{2}\left\Vert M-m\right\Vert \left[ \frac{1}{Q_{n}}%
\sum_{i=1}^{n}q_{i}\left\Vert z_{i}\right\Vert ^{2}-\left\Vert \frac{1}{Q_{n}%
}\sum_{i=1}^{n}q_{i}z_{i}\right\Vert ^{2}\right] ^{\frac{1}{2}}  \notag \\
& \leq \frac{1}{4}\left\Vert M-m\right\Vert \left\Vert Z-z\right\Vert . 
\notag
\end{align}
\end{corollary}

\begin{remark}
Note that the inequality between the first term and the last term in (\ref%
{3.9}) was firstly proved in \cite[Theorem 4.1]{D2}. Consequently, the above
corollary provides an improvement of the reverse of Jensen's inequality
established in \cite{D2}.
\end{remark}

\end{document}